\documentclass[a4paper,11pt]{amsart}
\usepackage{amsthm,amssymb,latexsym, mathrsfs}
\input amssym.def

\title[Two-weight norm estimates for maximal function ]{On two-weight norm estimates for multilinear fractional maximal function}
\newtheorem{theorem}{T{\hskip 0pt\footnotesize\bf HEOREM}}[section]
\newtheorem{lemma}[theorem]{L{\hskip 0pt\footnotesize\bf EMMA}}
\newtheorem{proposition}[theorem]{P{\hskip 0pt\footnotesize\bf ROPOSITION}}
\newtheorem{definition}[theorem]{D{\hskip 0pt\footnotesize\bf EFINITION}}
\newtheorem{corollary}[theorem]{C{\hskip 0pt\footnotesize\bf OROLLARY}}
\newtheorem{remark}[theorem]{R{\hskip 0pt\footnotesize\bf EMARK}}

\newcommand{\Proof}{\noindent{\bf P{\footnotesize \small\bf roof}:}}

\newcommand{\supp}{\mathrm{supp}}

\newcommand{\bprop} {\begin{proposition}}
\newcommand{\eprop} {\end{proposition}}
\newcommand{\btheo} {\begin{theorem}}
\newcommand{\etheo} {\end{theorem}}
\newcommand{\blem} {\begin{lemma}}
\newcommand{\elem} {\end{lemma}}
\newcommand{\bcor} {\begin{corollary}}
\newcommand{\ecor} {\end{corollary}}

\newcommand{\Be}{\begin{equation}}
\newcommand{\Ee}{\end{equation}}
\newcommand{\Bea}{\begin{eqnarray}}
\newcommand{\Eea}{\end{eqnarray}}
\newcommand{\Bes}{\begin{equation*}}
\newcommand{\Ees}{\end{equation*}}
\newcommand{\Beas}{\begin{eqnarray*}}
\newcommand{\Eeas}{\end{eqnarray*}}
\newcommand{\Ba}{\begin{array}}
\newcommand{\Ea}{\end{array}}

 \scrollmode

\begin{document}

\author[Beno\^it F. Sehba]{Beno\^it F. Sehba}
\address{Beno\^it Florent Sehba,  Department of Mathematics\\ University of Ghana, P.O.Box 62\\ Legon, Accra, Ghana}
\email{bfsehba@ug.edu.gh}
\keywords{Carleson embeddings, Fractional maximal function, $A_p$ weight}
\subjclass[2000]{Primary: 42B25 , Secondary: 42B20, 42B35 }


\begin{abstract}
We prove some Sawyer-type characterizations for multilinear fractional maximal function for the upper triangle case. We also provide some two-weight norm estimates for this operator. As one of the main tools, we use an extension of the usual Carleson Embedding that is an analogue of the P. L. Duren extension of the Carleson Embedding for measures.
\end{abstract}
\maketitle
\section{Introduction}
All over the text, $\mathbb R^n$ will be the $n$-dimensional real Euclidean space; all the cubes considered are non-degenerate with sides parallel to the coordinate axes. If $Q$ is a cube, then we denote by $|Q|$ its Lebesgue measure. When $\omega$ is a weight on $\mathbb R^n$, we write $\omega(Q):=\int_Q\omega(x)dx$. Given an exponent $1< p<\infty$, we denote by $p'$ its conjugate; that is $pp'=p+p'$.

An important question in modern harmonic analysis is given an operator $T$  determine the pairs of weights $(\omega,\sigma)$ such that
\Be\label{eq:form1}
\|Tf\|_{L^p(\omega)}\le C(\omega,\sigma)\|f\|_{L^p(\sigma)}
\Ee
or more generally,
\Be\label{eq:form2}
\|Tf\|_{L^q(\omega)}\le C(\omega,\sigma)\|f\|_{L^p(\sigma)}.
\Ee
When $T$ is the Hardy-Littlewood maximal function $M$, a complete answer was provided in the case $\omega=\sigma$ by B. Muckenhoupt \cite{Muckenhoupt}  who proved that (\ref{eq:form1}) holds for $M$ if and only if $\sigma$ satisfies the so-called $A_p$ condition. That is
\Be\label{eq:Apcondi}
[\sigma]_{A_p}:=\sup_{Q}\left(\frac{1}{|Q|}\int_Q\sigma\right)\left(\frac{1}{|Q|}\int_Q\sigma^{1-p'}\right)^{p-1}<\infty.
\Ee
Note however that this question in its general form is difficult.

Recall that the fractional maximal function is defined by
\Be\label{eq:fractmaxfunct}
M_{\alpha}f(x):=\sup_{Q}\frac{\chi_Q(x)}{|Q|^{1-\alpha /n}}\int_Q|f(y)|dy
\Ee
provided $0\le \alpha<n$. When $\alpha=0$, this is just the Hardy-Littlewood maximal function. B. Muckenhoupt and R. Wheeden \cite{MuckenhouptWhee} proved that for $1<p<\frac{n}{\alpha}$, and $\frac{1}{q}=\frac{1}{p}-\frac{\alpha}{n}$, $M_\alpha$ is bounded from $L^p(\sigma^p)$ to $L^q(\sigma^q)$ if and only if $$[\sigma]_{A_{p,q}}:=\sup_Q\left(\frac{1}{|Q|}\int_Q\sigma^q\right)^{1/q}\left(\frac{1}{|Q|}\int_Q\sigma^{-p'}\right)^{1/{p'}}<\infty.$$
In \cite{Sawyer}, E. Sawyer provided a general criterium for the maximal function to be bounded from $L^p(\sigma)$ to $L^q(\omega)$ which can be rephrased as follows (see \cite{Moen1} for the norm estimate).
\btheo\label{thm:sawyer}
Let $\sigma$ and $\omega$ be  two weights on $\mathbb R^n$, and $1< p\le q<\infty$. Then the following are equivalent
\begin{itemize}
\item[(i)] There exists a constant $C_1>0$ such that $$\int_{\mathbb R^n}\left(M_\alpha(\sigma f)(x)\right)^q\omega(x)dx\le C_1\|f\|_{p,\sigma}^q.$$
\item[(ii)] There exists a constant $C_2>0$ such that for any cube $Q$, $$\int_Q\left(M_\alpha(\sigma \chi_Q)(x)\right)^q\omega(x)dx\le C_2\left(\sigma(Q)\right)^{q/p}.$$
\end{itemize}
Moreover, if $$[\sigma,\omega]_{S_{p,q,\alpha}}:=\sup_{Q}\left(\frac{\int_Q\left(M_\alpha(\sigma \chi_Q)(x)\right)^q\omega(x)dx}{\left(\sigma(Q)\right)^{q/p}}\right)^{1/q},$$ then
\Be\label{eq:Spqestim}\|M\|_{L^p(\sigma)\rightarrow L^q(\omega)}\backsimeq C(\alpha,n)\left(\frac{q}{p}\right)^{1/q}p'[\sigma,\omega]_{S_{p,q,\alpha}}.\Ee
\etheo

We are interested in the multilinear analogues of the maximal operators above. For $m$ a given positive integer, the multilinear fractional maximal function is defined by
$$\mathcal {M}_\alpha \vec {f}(x):=\sup_{Q}|Q|^{\alpha /n}\prod_{i=1}^m\frac{\chi_Q(x)}{|Q|}\int_Q|f_i(y)|dy$$
provided $0\le \alpha<mn$. Here $\vec f=(f_1,\cdots,f_m)$ where the $f_i$s are measurable functions. When $\alpha=0$, $\mathcal M_0=\mathcal M$ is the multilinear Hardy-Littlewood maximal function. Note that these operators are related to  multilinear Calder\'on-Zygmund theory and the study of multilinear fractional integral operators \cite{Grafakos, GrafaKalton, GrafaTorres1, GrafaTorres2, KenigStein, LerOmbroPerezetal, Moen2}.

One of the first interest in this setting in relation with the linear case results commented above was the extension of Muckenhoupt result for $M$ and the Muckenhoupt-Wheeden result for $M_\alpha$. In \cite{LerOmbroPerezetal}, the authors provided the right characterization of weights for $\mathcal M$ to be bounded. More precisely, given $1<p_1<\cdots<p_m<\infty$, $\frac{1}{p}=\frac{1}{p_1}+\cdots+\frac{1}{p_m}$, and for $\vec {\sigma}=(\sigma_1,\cdots, \sigma_m)$ where the $\sigma_i$s are weights, define the weight $\omega_{\vec \sigma}$ by $\omega_{\vec \sigma}=\prod_{i=1}^m\sigma_i^{p/{p_i}}$. They proved that $\mathcal M$ is bounded from $L^{p_1}(\sigma_1)\times\cdots\times L^{p_m}(\sigma_m)$ to $L^p(\omega_{\vec \sigma})$ if and only if $$[\vec \sigma]_{A_{\vec P}}:=\sup_Q\left(\frac{1}{|Q|}\int_Q\omega_{\vec \sigma}\right)^{1/p}\prod_{i=1}^m\left(\frac{1}{|Q|}\int_Q\sigma_i^{1-p_i'}\right)^{1/{p_i'}}<\infty.$$
Building over these new lights, K. Moen considered the same question for $\mathcal M_\alpha$ in \cite{Moen2}. He proved that $\mathcal M_\alpha$ is bounded from $L^{p_1}(\sigma_1^{p_1})\times\cdots\times L^{p_m}(\sigma_m^{p_m})$ to $L^q(\omega_{\vec \sigma}^q)$ with $\frac{1}{q}=\frac{1}{p}-\frac{\alpha}{n}$ if and only if
$$[\vec \sigma]_{A_{\vec {P},q}}:=\sup_Q\left(\frac{1}{|Q|}\int_Q\left(\prod_{i=1}^m\sigma_i\right)^q\right)^{1/q}\prod_{i=1}^m\left(\frac{1}{|Q|}\int_Q\sigma_i^{-p_i'}\right)^{1/{p_i'}}<\infty.$$

The other question that comes in mind is the extension of Sawyer result to this setting. In \cite{ChenDamian} and \cite{LiXueetal}, the authors dealt with this question for $\alpha=0$ but under the assumption that the weights satisfy a kind of reverse H\"older inequality and monotone property respectively. K. Li and W. Sun \cite{LiSun} managed to extend the Sawyer characterization for the boundedness of $\mathcal M_\alpha$ from $L^{p_1}(\sigma_1)\times\cdots\times L^{p_m}(\sigma_m)$ to $L^q(\omega)$ for $\max\{p_1,\cdots,p_m\}\le q<\infty$. They proved the following.
\btheo\label{thm:LiSun} Given a nonnegative integer $m$, and $1<p_1,\cdots,p_m<\infty$; suppose that $0\le \alpha<mn$, $\frac{1}{p}=\frac{1}{p_1}+\cdots+\frac{1}{p_m}$ and $\max\{p_1,\cdots,p_m\}\le q<\infty$. Let $\omega_1,\cdots,\omega_m$ and $v$ be weights and put $\sigma_i=\omega_i^{1-p_i'}$, $i=1,\cdots,m$. Define
$$[\vec {\omega},v]_{S_{\vec {P},q}}:=\sup_{Q}\frac{\left(\int_Q\left(\mathcal {M}_\alpha(\sigma_1\chi_Q,\cdots,\sigma_m\chi_Q)(x)\right)^qv(x)dx\right)^{1/q}}{\prod_{i=1}^m\sigma_i(Q)^{1/{p_i}}}.$$
Then $\mathcal M_\alpha$ is bounded from $L^{p_1}(\omega_1)\times\cdots\times L^{p_m}(\omega_m)$ to $L^q(v)$ if and only if $[\vec {\omega},v]_{S_{\vec {P},q}}$ is finite. Moreover,
$$\|\mathcal M_\alpha\|_{\left(\prod_{i=1}^mL^{p_i}(\omega_i)\right)\rightarrow L^q(v)}\backsimeq [\vec {\omega},v]_{S_{\vec {P},q}}.$$

\etheo

The condition $[\vec {\omega},v]_{S_{\vec {P},q}}<\infty$ is necessary in general but the assumption $q\ge \max\{p_1,\cdots,p_m\}$ makes the result above restrictive. One might be interested in knowing if it is possible to remove this assumption and may be replace it by $q\ge p$, with $\frac{1}{p}=\frac{1}{p_1}+\cdots+\frac{1}{p_m}$. Before going ahead on this question, let us first observe the following general sufficient condition.
\bprop\label{prop:suffsawyertype} Given a nonnegative integer $m$, $1<p_1,\cdots,p_m<\infty$. Suppose that $0\le \alpha<mn$, $\frac{1}{p}=\frac{1}{p_1}+\cdots+\frac{1}{p_m}$ and $p\le q<\infty$. Let $\omega_1,\cdots,\omega_m$ and $v$ be weights and put $\sigma_i=\omega_i^{1-p_i'}$, $i=1,\cdots,m$ and $\nu_{\vec \omega}=\prod_{i=1}^m\sigma_i^{\frac{p}{p_i}}$. Define
$$[\nu_{\vec \omega},v]_{S_{\vec {P},q}}:=\sup_{Q}\frac{\left(\int_Q\left(\mathcal {M}_\alpha(\sigma_1\chi_Q,\cdots,\sigma_m\chi_Q)(x)\right)^qv(x)dx\right)^{1/q}}{\left(\nu_{\vec \omega}(Q)\right)^{1/p}}.$$
Then $\mathcal M_\alpha$ is bounded from $L^{p_1}(\omega_1)\times\cdots\times L^{p_m}(\omega_m)$ to $L^q(v)$  if $[\nu_{\vec \omega},v]_{S_{\vec {P},q}}$ is finite. Moreover,
$$\|\mathcal M_\alpha\|_{\left(\prod_{i=1}^mL^{p_i}(\omega_i)\right)\rightarrow L^q(v)}\lesssim [\nu_{\vec \omega},v]_{S_{\vec {P},q}}.$$

\eprop
It comes that if the weight $\omega_i$, $i=1,2,\cdots,m$ are such that for any cube $Q$,
\Be\label{eq:RH}
\prod_{i=1}^m\sigma_i(Q)^{p/{p_i}}=\prod_{i=1}^m\left(\int_Q\sigma_i(x)dx\right)^{p/p_i}\lesssim \int_Q\left(\prod_{i=1}^m\sigma_i^{\frac{p}{p_i}}\right)(x)dx=\nu_{\vec \omega}(Q),
\Ee
then the equivalence $[\nu_{\vec \omega},v]_{S_{\vec {P},q}}\backsimeq [\vec {\omega},v]_{S_{\vec {P},q}}$ holds and consequently, Theorem \ref{thm:LiSun}
holds without the restriction $q\ge \max\{p_1,\cdots,p_m\}$ but with this time $p\le q$. That is the following holds.
\btheo\label{thm:ChenDamianextended} Given a nonnegative integer $m$, $1<p_1,\cdots,p_m<\infty$. Suppose that $0\le \alpha<mn$, $\frac{1}{p}=\frac{1}{p_1}+\cdots+\frac{1}{p_m}$ and $p\le q<\infty$. Let $\omega_1,\cdots,\omega_m$ and $v$ be weights and put $\sigma_i=\omega_i^{1-p_i'}$, $i=1,\cdots,m$. Suppose that the weights $\sigma_i$, $i=1,\cdots,m$ are such that (\ref{eq:RH}) holds and define
$$[\nu_{\vec \omega},v]_{S_{\vec {P},q}}=\sup_{Q}\frac{\left(\int_Q\left(\mathcal {M}_\alpha(\sigma_1\chi_Q,\cdots,\sigma_m\chi_Q)(x)\right)^qv(x)dx\right)^{1/q}}{\left(\nu_{\vec \omega}(Q)\right)^{1/p}}.$$
Then $\mathcal M_\alpha$ is bounded from $L^{p_1}(\omega_1)\times\cdots\times L^{p_m}(\omega_m)$ to $L^q(v)$  if and only if $[\nu_{\vec \omega},v]_{S_{\vec {P},q}}$ is finite. Moreover,
$$\|\mathcal M_\alpha\|_{\left(\prod_{i=1}^mL^{p_i}(\omega_i)\right)\rightarrow L^q(v)}\backsimeq [\nu_{\vec \omega},v]_{S_{\vec {P},q}}.$$
\etheo
\vskip .1cm
Condition (\ref{eq:RH}) was used and named reverse H\"older inequality $RH_{\vec P}$ in \cite{ChenDamian, ChenLiu}. In \cite{ChenDamian}, the authors obtained Theorem \ref{thm:ChenDamianextended} for $\alpha=0$ and $p=q$, but it is hard to provide examples of family of weights for which (\ref{eq:RH}) holds. Nevertheless one can check that for
$\sigma_1=\sigma_2=\cdots=\sigma_m=\sigma$, we have the inequality (\ref{eq:RH}) and in this case, the following result.

\bcor\label{thm:Sehba}
Given a nonnegative integer $m$, $1<p_1,\cdots,p_m<\infty$. Suppose that $0\le \alpha<mn$, $\frac{1}{p}=\frac{1}{p_1}+\cdots+\frac{1}{p_m}$ and $p\le q<\infty$. Let $\sigma$ and $\omega$ be weights. Define
$$[\sigma,\omega]_{ {S}_{\vec {P},q}}=\sup_{Q}\frac{\left(\int_Q\left(\mathcal {M}_\alpha(\sigma\chi_Q,\cdots,\sigma\chi_Q)(x)\right)^q\omega(x)dx\right)^{1/q}}{\sigma(Q)^{1/{p}}}.$$
Then $\mathcal M_\alpha$ is bounded from $L^{p_1}(\sigma^{-p_1/p_1'})\times\cdots\times L^{p_m}(\sigma^{-p_m/p_m'})$ to $L^q(\omega)$ if and only if $[\sigma,\omega]_{\mathcal S_{\vec {P},q}}$ is finite. Moreover,
$$\|\mathcal M_\alpha\|_{\left(\prod_{i=1}^mL^{p_i}(\sigma^{-p_i/p_i'})\right)\rightarrow L^q(\omega)}\backsimeq [\omega,v]_{S_{\vec {P},q}}.$$

\ecor
Recall that the $\mathcal A_\infty$ class of Hru\'s\'cev (\cite{Hurscev}) consists of weights $\omega$ satisfying
\Be\label{AinftyHurscev0}
[\omega]_{\mathcal A_\infty}:=\sup_{Q}\left(\frac{1}{|Q|}\int_Q\omega\right)\exp\left(\frac{1}{|Q|}\int_Q\log \omega^{-1}\right)<\infty.
\Ee
If is easy to check that for $\sigma_1,\cdots,\sigma_m\in\mathcal {A}_{\infty} $, and for any cube $Q$, 
$$\prod_{i=1}^m\sigma_i(Q)^{p/{p_i}}\lesssim \left(\prod_{i=1}^m[\sigma_i]_{\mathcal A_\infty}\right)\int_Q\left(\prod_{i=1}^m\sigma_i^{\frac{p}{p_i}}\right)(x)dx $$
(see \cite{WangYi}). It follows that we also have the following result.
\bcor\label{thm:ChenDamianextended1} Given a nonnegative integer $m$, $1<p_1,\cdots,p_m<\infty$. Suppose that $0\le \alpha<mn$, $\frac{1}{p}=\frac{1}{p_1}+\cdots+\frac{1}{p_m}$ and $p\le q<\infty$. Let $\omega_1,\cdots,\omega_m$ and $v$ be weights and put $\sigma_i=\omega_i^{1-p_i'}$, $i=1,\cdots,m$. Suppose that the weights $\sigma_i$, $i=1,\cdots,m$ are  in the class $\mathcal A_\infty$ and define
$$[\vec {\omega},v]_{S_{\vec {P},q}}=\sup_{Q}\frac{\left(\int_Q\left(\mathcal {M}_\alpha(\sigma_1\chi_Q,\cdots,\sigma_m\chi_Q)(x)\right)^qv(x)dx\right)^{1/q}}{\prod_{i=1}^m \sigma_i(Q)^{1/p_i}}.$$
Then $\mathcal M_\alpha$ is bounded from $L^{p_1}(\omega_1)\times\cdots\times L^{p_m}(\omega_m)$ to $L^q(v)$  if and only if $[\vec {\omega},v]_{S_{\vec {P},q}}$ is finite. Moreover,
$$\|\mathcal M_\alpha\|_{\left(\prod_{i=1}^mL^{p_i}(\omega_i)\right)\rightarrow L^q(v)}\backsimeq [\vec {\omega},v]_{S_{\vec {P},q}}.$$
\ecor
%

To prove Proposition \ref{prop:suffsawyertype}, one first need to observe that the matter can be reduced to the associated dyadic maximal function. We then use an approach that can be traced back to \cite{Sawyer} and has been simplified in \cite{Cruz}, it consists in discritizing the integral $\int_{\mathbb R^d}\left(\mathcal M_{d,\alpha}(f)\right)^q\omega$ where $\mathcal M_{d,\alpha}$ stands for the multilinear dyadic fractional maximal function, using appropriate level sets and their decomposition into disjoint dyadic cubes. In the linear case (i.e. when $m=1$), one then uses an interpolation approach to get the embedding (see \cite{Cruz, LiSun}). This method still works in the multilinear case under further restrictions on the exponents that allow one to reduce the matter to a linear case and this is what happens exactly in the proof of Theorem \ref{thm:LiSun} in \cite{LiSun}. It is not clear how this can be done in general in the multilinear setting for the upper triangle case ($p< q$). To overcome this difficulty, we just extend the techniques used for the diagonal case ($p=q$) which reduce the matter to a Carleson embedding (see \cite{ChenDamian}). More precisely, we use the following extension of the usual Carleson embedding and its multilinear analogue.

\btheo\label{thm:Carlembed0}
Let $\sigma$ be a weight on $\mathbb R^n$ and $\alpha\ge 1$. Assume  $\{\lambda_Q\}_{Q\in \mathcal D}$ is a sequence of positive numbers indexed over the set of dyadic cubes $\mathcal D$ in $\mathbb R^n$. Then the following are equivalent.
\begin{itemize}
\item[(i)] There exists some constant $A>0$ such that for any cube $R\in \mathcal D$,
    \begin{equation*}\sum_{Q\subseteq R, Q\in \mathcal D}\lambda_Q\le A(\sigma(R))^\alpha.\end{equation*}
\item[(ii)] There exists a constant $B>0$ such that for all $p\in [1,\infty)$,
\begin{equation*}\sum_{Q\in \mathcal D}\lambda_Q|m_\sigma(f,Q)|^{p\alpha}\le B \|f\|_{p,\sigma}^{p\alpha}.\end{equation*}
\end{itemize}
\etheo

The above theorem for $\alpha=1$ is the usual Carleson Embedding Theorem (see \cite{HytPerez, Pereyra}). The case  $\alpha>1$ seems to be new and can be viewed as an analogue of P. L. Duren extension of Carleson embedding theorem for measures \cite{Duren}.The proof of Theorem \ref{thm:LiSun} is also simplified when combining the main idea of \cite{LiSun} and the extension of the Carleson embedding. For the proof of Prosition \ref{prop:suffsawyertype}, we will use a multilinear analogue of the above embedding.

\vskip .2cm
Our other interest in this paper is to provide sufficient conditions for $\mathcal M_\alpha$ to be bounded from $L^{p_1}(\sigma_1)\times\cdots\times L^{p_m}(\sigma_m)$ to $L^q(\omega)$. Usually, one expects conditions that have a form close to the $A_p$ characteristic of Muckenhoupt. This question is quite interesting in this research area as it is related to the same type of questions for singular operators and some questions arising from PDEs (see \cite{CruzMoen, GrafaTorres2,Moen2, MuckenhouptWhee,Perez} and the references therein).  Before going ahead on this question, we need more definitions and notations.

Given two weights $\omega$ and $\sigma$, we say they satisfy the joint $A_p$ condition for $1<p<\infty$  if \Be\label{eq:Apclass}[\omega,\sigma]_{A_p}:=\sup_{Q}\frac{\omega(Q)\sigma(Q)^{p-1}}{|Q|^p}<\infty.\Ee
Note that when $\sigma=\omega^{-1/(p-1)}$, this is just the definition of the $A_p$ class of Muckenhoupt. A new class of weights was recently introduced by T. Hyt\"onen and C. P\'erez \cite{HytPerez} and consists of pair of weights satisfying the condition
  \Be\label{eq:Bpclass}[\omega,\sigma]_{B_p}:=\sup_{Q}\frac{\omega(Q)\sigma(Q)^{p}}{|Q|^{p+1}}\exp\left(\frac{1}{|Q|}\int_Q\log \sigma^{-1}\right)<\infty.\Ee
We recall the definition of the $A_\infty$ class of Fujii-Wilson (\cite{Fujii, Lerner,wilson1,wilson2,wilson3}). We say a weight $\sigma$ belongs to $A_\infty$ if
\Be\label{eq:Ainftyclass}[\sigma]_{A_\infty}:=\sup_{Q}\frac{1}{\sigma(Q)}\int_Q M(\sigma \chi_Q)<\infty.\Ee
S. Buckley \cite{Buckley} obtained the following estimate for the maximal operator $M$: \Be\label{eq:Buckleyestim}\|M\|_{L^p(\sigma)\rightarrow L^p(\sigma)}\le Cp'[\sigma]_{A_p}^{1/(p-1)}.\Ee
This was recently improved by T. Hyt\"onen and C. P\'erez \cite{HytPerez} as follows.
\btheo\label{thm:hytperez}
Let $1<p<\infty$, and $\sigma$, $\omega$ two weights. Then
\Be\label{eq:hytperez1}
\|M(f\sigma)\|_{p,\omega}\le C(\alpha,n)p'\left([\sigma,\omega]_{B_p}\right)^{1/p}\|f\|_{p,\sigma}
\Ee
and
\Be\label{eq:hytperez2}
\|M(f\sigma)\|_{p,\omega}\le C(\alpha,n)p'\left([\sigma,\omega]_{A_p}[\sigma]_{A_\infty}\right)^{1/p}\|f\|_{p,\sigma}.
\Ee
\etheo
Estimate (\ref{eq:hytperez2}) in the case $p=2$ is actually attributed to A. K. Lerner and S. Ambrosi \cite{LerOmbro}. To find the corresponding estimates in the $L^p-L^q$ case, we need to introduce adapted classes of weights that generalize the above ones. For $1<p_1,\cdots,p_m,q<\infty$, $\vec {P}=(p_1,\cdots,p_m)$, we say the weights $\vec {\sigma}=(\sigma_1,\cdots,\sigma_m)$ and $\omega$ satisfy the joint conditions $A_{\vec P,q}$ and $B_{\vec P,q}$ if
\Be\label{eq:Apqclass}[\vec {\sigma},\omega]_{A_{\vec P,q}}:=\sup_{Q}\frac{\omega(Q)^{p/q}\prod_{i=1}^m\sigma_i(Q)^{p/{p_i'}}}{|Q|^{p(m-\alpha /n)}}<\infty
\Ee
and
\Be\label{eq:Bpqclass}[\vec {\sigma},\omega]_{B_{\vec P,q}}:=\sup_{Q}\frac{\omega(Q)^{p/q}\prod_{i=1}^m\sigma_i(Q)^{p}}{|Q|^{p(m-\alpha /n)+1}}\prod_{i=1}^m\left(\exp\left(\frac{1}{|Q|}\int_Q\log \sigma^{-1}\right)\right)^{p/p_i}<\infty.\Ee

One easily checks the following inequalities
$$[\vec {\sigma},\omega]_{A_{\vec P,q}}\le [\vec {\sigma},\omega]_{B_{\vec P,q}}\le [\vec {\sigma},\omega]_{A_{\vec P,q}}\prod_{i=1}^m[\sigma_i]_{\mathcal A_\infty}^{p/p_i}.$$
Let us also introduce the multilinear $A_\infty$ class of Chen-Damian \cite{ChenDamian}. That is the class of weights $\vec {\sigma}=(\sigma_1,\cdots,\sigma_m)$ such that
$$[\vec {\sigma}]_{W_{\vec P}^\infty}:=\sup_{Q}\frac{\int_Q\prod_{i=1}^mM(\sigma_i\chi_Q)(x)^{p/p_i}dx}{\int_Q\prod_{i=1}^m\sigma_i^{p/p_i}(x)dx}<\infty,$$
where $\frac{1}{p}=\frac{1}{p_1}+\cdots+\frac{1}{p_m}$.
Our corresponding result is the following.
\btheo\label{thm:main}
Let $1<p_1,\cdots,p_m, q<\infty$, and $\vec {\sigma}=(\sigma_1,\cdots,\sigma_m)$, $\omega$ be weights. Put
$\frac{1}{p}=\frac{1}{p_1}+\cdots+\frac{1}{p_m}$ and assume that $p\le q$.
Then
\Be\label{eq:main1}
\|\mathcal M_\alpha(\sigma_1f_1,\cdots,\sigma_mf_m)\|_{q,\omega}\le C(n,p,q)\left([\vec {\sigma},\omega]_{B_{\vec P,q}}\right)^{1/p}\prod_{i=1}^m\|f_i\|_{p_i,\sigma_i},
\Ee

\Be\label{eq:main2}
\|\mathcal M_\alpha(\sigma_1f_1,\cdots,\sigma_mf_m)\|_{q,\omega}\le C(n,p,q)\left([\vec {\sigma},\omega]_{A_{\vec P,q}}\right)^{1/p}\left(\prod_{i=1}^m[\sigma_i]_{A_\infty}^{1/p_i}\right)\prod_{i=1}^m\|f_i\|_{p_i,\sigma_i}
\Ee
and
\Be\label{eq:main3}
\|\mathcal M_\alpha(\sigma_1f_1,\cdots,\sigma_mf_m)\|_{q,\omega}\le C(n,p,q)\left([\vec {\sigma},\omega]_{A_{\vec P,q}}[\vec {\sigma}]_{W_{\vec P}^\infty}\right)^{1/p}\prod_{i=1}^m\|f_i\|_{p_i,\sigma_i}.
\Ee
\etheo
We observe that when $\alpha=0$ and $p=q$, that is for the multilinear Hardy-Littlewood maximal function, inequalities (\ref{eq:main1}) and (\ref{eq:main3}) were proved in \cite{ChenDamian}, while (\ref{eq:main2}) was obtained in \cite{DamLerPer}. Sharp norm estimates of the fractional maximal function are considered in \cite{LiMoenetal}, these estimates are similar to (\ref{eq:main2}) with a modification of the power on $[\sigma_i]_{A_\infty}$. An extension of the Buckley estimate (\ref{eq:Buckleyestim}) to the multilinear maximal function is given in \cite{DamLerPer}.
\vskip .2cm

To prove Theorem \ref{thm:main}, one first needs to observe as above that one only needs to consider the case of the dyadic maximal function. Next to estimate the norm of the dyadic maximal function, we proceed essentially as for Proposition \ref{prop:suffsawyertype}. For some other sufficient conditions of this type, we refer the reader to the following and the references therein \cite{ChenDamian,LerOmbroPerezetal,LiMoenetal,Moen1,Moen2}.
\vskip .2cm
The paper is organized as follows, in the next section, we introduce an extension of the usual notion of Carleson sequences, and provide equivalent characterizations. In section 3, we prove Proposition \ref{prop:suffsawyertype} and simplify the proof of Theorem \ref{thm:LiSun}. Theorem \ref{thm:main} is proved in the last section. Some steps in our proofs are known by the specialists but we write them down so that the reader can easily follow us.
\vskip .2cm
All  over the text, $C$ will denote a constant not necessarily the same at each occurrence. We write $C(\alpha,n,\cdots)$ to emphasize on the fact that our constant depends on the parameters $\alpha,n,\cdots$. As usual, given two positive quantities $A$ and $B$, the notation $A\lesssim B$ (resp. $B\lesssim A$) will mean that there is an universal constant $C>0$ such that $A\le CB$ (resp. $B\le CA$). When $A\lesssim B$ and $B\lesssim A$, we write $A\backsimeq B$ and say $A$ and $B$ are equivalent.
\section{$(\alpha,\sigma)$-Carleson sequences}
We introduce a more general notion of Carleson sequences, provide equivalent definitions and applications.
\subsection{Definitions and results}
We have the following general definition of Carleson sequences.
\begin{definition}\label{def:alphacarlseq}
Given a weight $\sigma$ and a number $\alpha\ge 1$, a sequence of positive numbers $\{\lambda_Q\}_{Q\in \mathcal D}$ indexed over the set of dyadic cubes $\mathcal D$ in $\mathbb R^n$ is called a $(\alpha,\sigma)$-Carleson sequence if there exists a constant $A>0$ such that for any cube $R\in \mathcal D$,
\Be
\sum_{Q\subseteq R}\lambda_Q\le A(\sigma(R))^\alpha.
\Ee
\end{definition}
We call Carleson constant of the sequence $\{\lambda_Q\}_{Q\in \mathcal D}$, the smallest constant in the above definition and denote it by $A_{Carl}$ when there is no ambiguity. When $\sigma\equiv 1$ and $\alpha\ge 1$, we speak of $\alpha$-Carleson sequences. In particular when $\alpha=1$, we just call them Carleson sequences as usual.

Let us introduce some notations. For $f\in L^p(\omega)$, $$\|f\|_{p,\omega}^p:=\int_{\mathbb R^n}|f(x)|^p\omega(x)dx$$ and $$m_\omega(f,Q):=\frac{1}{\omega(Q)}\int_Qf(x)\omega(x)dx$$ where $\omega(Q)=\int_Q\omega$. When $\omega\equiv 1$, we write
$m_Qf=m(f,Q)=m_\omega(f,Q)$.

We have the following characterization of $(\alpha,\sigma)$-Carleson sequences.
\btheo\label{thm:Carlembed1}
Let $\sigma$ be a weight on $\mathbb R^n$ and $\alpha\ge 1$. Assume  $\{\lambda_Q\}_{Q\in \mathcal D}$ is a sequence of positive numbers indexed over the set of dyadic cubes $\mathcal D$ in $\mathbb R^n$. Then the following are equivalent.
\begin{itemize}
\item[(i)] $\{\lambda_Q\}_{Q\in \mathcal D}$ is a $(\alpha,\sigma)$-Carleson sequence, that is for some constant $A>0$ and for any cube $R\in \mathcal D$,
    \Be\label{eq:Carlcond1}\sum_{Q\subseteq R, Q\in \mathcal D}\lambda_Q\le A(\sigma(R))^\alpha.\Ee
\item[(ii)]There exists a constant $B>0$ such that for all $p\in [1,\infty)$,
\Be\label{eq:Carlembed1}\sum_{Q\in \mathcal D}\lambda_Q|m_\sigma(f,Q)|^{p\alpha}\le B \|f\|_{p,\sigma}^{p\alpha}.\Ee
\end{itemize}
\etheo
\begin{proof}
Let us recall that the dyadic Hardy-Littlewood maximal function with respect to the measure $\sigma$ is defined by
$$M_{d}^{\sigma}f:=\sup_{Q\in \mathcal D}\frac{\chi_Q}{\sigma(Q)}\int_Q|f|\sigma.$$
When $\sigma$ is the Lebesgue measure, we write $M_{d}^{\sigma}=M_d$. We recall the estimate
\Be\label{eq:dyamaxfunctestim}
\|M_{d}^{\sigma}f\|_{p,\sigma}\le p'\|f\|_{p,\sigma}.
\Ee
We will also need the following inequality.
\Be\label{eq:weakineq}
\lambda^p\sigma\left(\{x:M_{d}^{\sigma}f(x)>\lambda\}\right)\le \|M_{d}^{\sigma}f\|_{p,\sigma}^p.
\Ee

For the implication $(\textrm{ii})\Rightarrow (\textrm{i})$, we test $(\textrm{ii})$ with the function $f=\chi_R$ with $R\in \mathcal D$ to obtain
\Beas
\sum_{Q\subseteq R, Q\in \mathcal D}\lambda_Q &\le& \sum_{Q\in \mathcal D}\lambda_Q \left(m_\sigma(f,Q)\right)^{p\alpha}\\ &\le& B\|\chi_R\|_{p,\sigma}^{p\alpha}=B\left(\sigma(R)\right)^\alpha.
\Eeas
That is for any $R\in \mathcal D$
$$\sum_{Q\subseteq R, Q\in \mathcal D}\lambda_Q\le B\left(\sigma(R)\right)^\alpha$$
which is $(\textrm{i})$.

To prove that $(\textrm{i})\Rightarrow (\textrm{ii})$, it is enough by (\ref{eq:dyamaxfunctestim}) to prove the following.
\blem\label{lem:Carlembedlemma}
Let $\{\lambda_Q\}_{Q\in \mathcal D}$ and $\alpha\ge 1$.
Suppose that there exists a constant $A>0$ such that for any $R\in \mathcal D$,
$$\sum_{Q\subseteq R, Q\in \mathcal D}\lambda_Q\le A(\sigma(R))^\alpha.$$
Then for all $p\in [1,\infty)$,
\Be\label{eq:Carlembedmaxi1}\sum_{Q\in \mathcal D}\lambda_Q|m_\sigma(f,Q)|^{p\alpha}\le A\alpha \|M_{d}^{\sigma}f\|_{p,\sigma}^{p\alpha}.\Ee
\elem
\begin{proof}
We can suppose that $f>0$. As in the case of $\alpha=1$ in \cite{HytPerez}, we read $\sum_{Q\in \mathcal D}\lambda_Q\left(m_\sigma(f,Q)\right)^{p\alpha}$ as an integral over the measure space $(\mathcal D, \mu)$ built over the set of dyadic cubes $\mathcal D$, with $\mu$ the measure assigning to each cube $Q\in \mathcal D$ the measure $\lambda_Q$. Thus
\Beas
\sum_{Q\in \mathcal D}\lambda_Q\left(m_\sigma(f,Q)\right)^{p\alpha} &=& \int_0^\infty p\alpha t^{p\alpha-1}\mu\left(\{Q\in \mathcal D: m_\sigma(f,Q)>t\}\right)dt\\ &=& \int_0^\infty p\alpha t^{p\alpha-1}\mu(\mathcal D_t)dt,
\Eeas
$\mathcal D_t:=\{Q\in \mathcal D: m_\sigma(f,Q)>t\}$. Let $\mathcal {D}_t^*$ be the set of maximal dyadic cubes $R$ with respect to the inclusion so that $m_\sigma(f,R)>t$. Then $\mathcal {D}_t^*$ is the union of these disjoint maximal cubes and we have $\mathcal {D}_t^*=\{x\in \mathbb R^n: M_{d}^{\sigma}f(x)>t\}$. It follows from the hypothesis on the sequence $\{\lambda_Q\}_{Q\in \mathcal D}$ that
\Beas
\mu(\mathcal D_t) &=& \sum_{Q\in \mathcal D_t}\lambda_Q\le \sum_{R\in \mathcal {D}_t^*}\sum_{Q\subseteq R}\lambda_Q\\ &\le& A\sum_{R\in \mathcal {D}_t^*}\left(\sigma(R)\right)^\alpha\le A\left(\sum_{R\in \mathcal {D}_t^*}\sigma(R)\right)^\alpha\\ &\le& A\left(\sigma(\{M_{d}^{\sigma}f>t\})\right)^\alpha.
\Eeas
Hence using (\ref{eq:weakineq}), we obtain
\Beas
&&\sum_{Q\in \mathcal D}\lambda_Q|m_\sigma(f,Q)|^{p\alpha}\\ &\le& A\int_0^\infty p\alpha t^{p\alpha-1}\left(\sigma(\{M_{d}^{\sigma}f>t\})\right)^\alpha dt\\ &=& A\int_0^\infty p\alpha t^{p-1}\sigma(\{M_{d}^{\sigma}f>t\})\left(t^p\sigma(\{M_{d}^{\sigma}f>t\})\right)^{\alpha-1} dt\\ &\le& A\alpha\|M_{d}^{\sigma}f\|_{p,\sigma}^{p(\alpha-1)}\int_0^\infty p t^{p-1}\sigma(\{M_{d}^{\sigma}f>t\})dt\\ &\le& A\alpha\|M_{d}^{\sigma}f\|_{p,\sigma}^{p\alpha}.
\Eeas
\end{proof}
The proof is complete.
\end{proof}
The above theorem is clearly a generalization as taking $\alpha=1$ we get the well known Carleson embedding result (see \cite{HytPerez, Pereyra}). 

\begin{remark}
As a first application, we obtain a necessary and sufficient condition for the main paraproduct to be bounded from $L^p(\mathbb R)$ to $L^2(\mathbb R)$ for $1\le p\le 2$. Let us still denote by $\mathcal{D}$ the set of dyadic intervals in $\mathbb{R}$. Recall that given a dyadic interval $I$, the Haar function supported by $I$ is defined by 
$h_{I}(s)=|I|^{-\frac{1}{2}}(\chi_{I^+}(s)-\chi_{I^-}(s))$, where $I^-$ and $I^+$ are the left and the right halfs of $I$ respectively. For $\phi\in L^2(\mathbb R)$ with finite Haar expansion, the (main) paraproduct with symbol $\phi$ is the operator defined on $L^2(\mathbb R)$ by $$\Pi_\phi b(s):=\sum_{I\in \mathcal D}\langle \phi,h_I\rangle (m_Ib) h_I(s)$$
where $m_Ib=\frac{1}{|I|}\int_Ib(x)dx$. It is well known that the operator $\Pi_\phi$ is bounded on $L^p(\mathbb R)$ if and only the sequence $\{|\langle \phi,h_I\rangle|^2\}_{I\in \mathcal D}$ is a Carleson sequence. The following partial extension is a direct consequence of Theorem \ref{thm:Carlembed1}.
\bcor\label{cor:paraprod}
Let $\phi\in L^2(\mathbb R)$ and $1\le p\le 2$. Then $\Pi_\phi$ extends to a bounded operator from $L^p(\mathbb R)$ to $L^2(\mathbb R)$ if and only if

    \Be A:=\sup_{J\in \mathcal D}\frac{1}{|J|^{2/p}}\sum_{I\subseteq J, I\in \mathcal D}|\langle \phi,h_I\rangle|^2<\infty.\Ee

Moreover, $\|\Pi_\phi\|_{L^p(\mathbb R)\rightarrow L^2(\mathbb R)}\backsimeq A$.
\ecor
The higher dimensional version of the above corollary requires an adapted multivariable version of Theorem \ref{thm:Carlembed1}. This will be presented elsewhere.
\end{remark}

An alternative characterization of $(\alpha,\sigma)$-Carleson sequences is the following.
\btheo\label{thm:Carlembed2}
Let $N\ge 1$ be an integer,  $1\le p_j<\infty$, $j=1,\cdots,N$. Assume  that $\sigma$ is a weight on $\mathbb R^n$, and that there are $0<q_1,\cdots,q_N<\infty$ such that $\alpha=\sum_{j=1}^N\frac{q_j}{p_j}\ge1$. Then given a sequence $\{\lambda_Q\}_{Q\in \mathcal D}$  of positive numbers, the following are equivalent.
\begin{itemize}
\item[(i)] $\{\lambda_Q\}_{Q\in \mathcal D}$ is a $(\alpha,\sigma)$-Carleson sequence, that is for some constant $A>0$ and for any cube $R\in \mathcal D$,
    $$\sum_{Q\subseteq R, Q\in \mathcal D}\lambda_Q\le A(\sigma(R))^\alpha.$$
\item[(ii)] There exists a constant $B>0$ such that
\Be\label{eq:Carlembed2}\sum_{Q\in \mathcal D}\lambda_Q\prod_{j=1}^N|m_\sigma(f_j,Q)|^{q_j}\le B \prod_{j=1}^N\|f_j\|_{p_j,\sigma}^{q_j}.\Ee
\end{itemize}
\etheo
\begin{proof}
To prove that $(\textrm{ii})\Rightarrow (\textrm{i})$, take for $R\in \mathcal D$ given, $f_j=\chi_R$ for $j=1,\cdots, N$ and proceed as in the proof of Theorem \ref{thm:Carlembed1}. To prove that $(\textrm{i})\Rightarrow (\textrm{ii})$, it is enough to prove the following lemma which might be useful in some other circumstances.
\blem\label{lem:Carlembedlemma2}
Let $N\ge 1$ be an integer,  $1\le p_j<\infty$, $j=1,\cdots,N$. Assume  that $\sigma$ is a weight on $\mathbb R^n$, and that $0<q_1,\cdots,q_N<\infty$ so that $\alpha=\sum_{j=1}^N\frac{q_j}{p_j}\ge 1$. Then if $\{\lambda_Q\}_{Q\in \mathcal D}$  is a sequence of positive numbers such that there exists a constant $A>0$ so that for any cube $R\in \mathcal D$,
    $$\sum_{Q\subseteq R, Q\in \mathcal D}\lambda_Q\le A(\sigma(R))^\alpha,$$
then
\Be\label{eq:Carlembedmax2}\sum_{Q\in \mathcal D}\lambda_Q\prod_{j=1}^N|m_\sigma(f_j,Q)|^{q_j}\lesssim A\alpha \prod_{j=1}^N\|M_{d}^{\sigma}f_j\|_{p_j,\sigma}^{q_j}.\Ee
\elem
\begin{proof}
An application of H\"older's inequality and Lemma \ref{lem:Carlembedlemma} provide
\Beas
\sum_{Q\in \mathcal D}\lambda_Q\prod_{j=1}^N|m_\sigma(f_j,Q)|^{q_j} &\le& \prod_{j=1}^N\left(\sum_{Q\in \mathcal D}\lambda_Q|m_\sigma(f_j,Q)|^{p_j\alpha}\right)^{\frac{q_j}{\alpha p_j}}\\ &\le& \prod_{j=1}^N\left(A\alpha\|M_{d}^{\sigma}f_j\|_{p_j,\sigma}^{p_j\alpha}\right)^{\frac{q_j}{\alpha p_j}}\\ &\lesssim& A\alpha\prod_{j=1}^N\|M_{d}^{\sigma}f_j\|_{p_j,\sigma}^{q_j}.
\Eeas
\end{proof}
The proof is complete.
\end{proof}
\subsection{Another extension}

We end this section with the following extension of the multilinear Carleson embedding of \cite{ChenDamian}.
\blem\label{lem:Carlembedchendamian}
Let $N\ge 1$ be an integer,  $1\le p_i, \alpha<\infty$, $i=1,\cdots,N$. Assume  that $\sigma_i$, $i=1,\cdots, m$ are weights on $\mathbb R^n$, and put $\nu_{\vec \sigma}=\prod_{i=1}^m\sigma_i^{p/p_i}$ and $\frac{1}{p}=\frac{1}{p_1}+\cdots+\frac{1}{p_m}$. Then if $\{\lambda_Q\}_{Q\in \mathcal D}$  is a sequence of positive numbers such that there exists a constant $A>0$ so that for any cube $R\in \mathcal D$,
    $$\sum_{Q\subseteq R, Q\in \mathcal D}\lambda_Q\le A(\nu_{\vec \sigma}(R))^\alpha,$$
then
\begin{align}\label{eq:Carlembedchendamian}\sum_{Q\in \mathcal D}\lambda_Q\left|\prod_{i=1}^Nm_{\sigma_i}(f_i,Q)\right|^{p\alpha} &\le A\alpha\|M_{d}^{\vec \sigma}(\vec f)\|_{p,\nu_{\vec \sigma}}^{p\alpha}\\ &\lesssim A\alpha \prod_{i=1}^N\|M_d^{\sigma_i}f_i\|_{p_i,\sigma_i}^{p\alpha}\notag\\ &\lesssim A\alpha \prod_{i=1}^N\|f_i\|_{p_i,\sigma_i}^{p\alpha},\notag
\end{align}
where $$M_{d}^{\vec \sigma}(\vec f)=\sup_{Q\in \mathcal D}\prod_{i=1}^m\frac{\chi_Q}{\sigma_i(Q)}\int_Q|f_i|\sigma_i(x)dx,$$
$\vec f=(f_i,\cdots,f_m)$.

\elem
\begin{proof}
We read $\sum_{Q\in \mathcal D}\lambda_Q\left(\prod_{i=1}^N|m_{\sigma_i}(f_i,Q)|\right)^{p\alpha}$ as an integral over the measure space $(\mathcal D, \mu)$ built over the set of dyadic cubes $\mathcal D$, with $\mu$ the measure assigning to each cube $Q\in \mathcal D$ the measure $\lambda_Q$. Thus
\Beas
\sum_{Q\in \mathcal D}\lambda_Q\left(\prod_{i=1}^N|m_{\sigma_i}(f_i,Q)|\right)^{p\alpha} &=& \int_0^\infty p\alpha t^{p\alpha-1}\mu\left(\{Q\in \mathcal D: \prod_{i=1}^N|m_{\sigma_i}(f_i,Q)|>t\}\right)dt\\ &=& \int_0^\infty p\alpha t^{p\alpha-1}\mu(\mathcal D_t)dt,
\Eeas
$D_t:=\{Q\in \mathcal D: \prod_{i=1}^N|m_{\sigma_i}(f_i,Q)|>t\}$. Let $\mathcal {D}_t^*$ be the set of maximal dyadic cubes $R$ with respect to the inclusion with $\prod_{i=1}^m\frac{\chi_Q}{\sigma_i(Q)}\int_Q|f_i|\sigma_i(x)dx>t$. Then $\mathcal {D}_t^*$ is the union of these disjoint maximal cubes and we have $\mathcal {D}_t^*=\{x\in \mathbb R^n: M_{d}^{\vec {\sigma}}(\vec {f})(x)>t\}$. It follows from the hypothesis on the sequence $\{\lambda_Q\}_{Q\in \mathcal D}$ that
\Beas
\mu(\mathcal D_t) &=& \sum_{Q\in \mathcal D_t}\lambda_Q\le \sum_{R\in \mathcal {D}_t^*}\sum_{Q\subseteq R}\lambda_Q\\ &\le& A\sum_{R\in \mathcal {D}_t^*}\left(\nu_{\vec \sigma}(R)\right)^\alpha\le A\left(\sum_{R\in \mathcal {D}_t^*}\nu_{\vec \sigma}(R)\right)^\alpha\\ &\le& A\left(\nu_{\vec \sigma}(\{M_{d}^{\vec {\sigma}}f>t\})\right)^\alpha.
\Eeas
Hence using (\ref{eq:weakineq}) applied to $M_{d}^{\vec \sigma}(\vec f)$, we obtain
\Beas
\sum_{Q\in \mathcal D}\lambda_Q\left(\prod_{i=1}^N|m_{\sigma_i}(f_i,Q)|\right)^{p\alpha} &\le& A\int_0^\infty p\alpha t^{p\alpha-1}\left(\sigma(\{M_{d}^{\vec {\sigma}}(\vec f)>t\})\right)^\alpha dt\\ &=& A\int_0^\infty p\alpha t^{p-1}\sigma(\{M_{d}^{\vec {\sigma}}(\vec f)>t\})\left(t^p\sigma(\{M_{d}^{\vec {\sigma}}(\vec f)>t\})\right)^{\alpha-1} dt\\ &\le& A\alpha\|M_{d}^{\vec {\sigma}}(\vec f)\|_{p,\sigma}^{p(\alpha-1)}\int_0^\infty p t^{p-1}\sigma(\{M_{d}^{\vec {\sigma}}(\vec f)>t\})dt\\ &\le& A\alpha\|M_{d}^{\vec {\sigma}}(\vec f)\|_{p,\sigma}^{p\alpha}.
\Eeas
The second inequality in (\ref{eq:Carlembedchendamian}) follows from the H\"older inequality while the third follows from (\ref{eq:dyamaxfunctestim}). The proof is complete.
\end{proof}
\section{Sawyer-type two-weight characterization}
Let us consider the following family of dyadic grids in $\mathbb R^n$. $$\mathcal D^\beta:=\{2^{-k}\left([0,1)^n+m+(-1)^k\beta\right):k\in \mathbb Z, m\in \mathbb {Z}^m\};\,\,\,\beta\in \{0,\frac{1}{3}\}^n.$$
For $\beta={\bf 0}=(0,0,\cdots.0)$, we write $\mathcal D^{{\bf 0}}=\mathcal D$. The dyadic multilinear fractional maximal function with respect to the grid $\mathcal D^\beta$ is defined by
$$\mathcal M_{d,\alpha}^{\beta}f(x):=\sup_{Q\in \mathcal D^\beta}\prod_{i=1}^m\frac{\chi_{Q}(x)}{|Q|^{1-\frac{\alpha}{mn}}}\int_{Q_\beta}|f_i(y)|dy.$$
When $\alpha=0$, this is just the dyadic multilinear Hardy-Littlewood maximal function denoted here $\mathcal M_d^\beta$. When $\beta={\bf 0}$, we write $\mathcal M_{d,\alpha}$ and $\mathcal M_{d}$ for $\mathcal M_{d,\alpha}^\beta$ and $\mathcal M_{d}^\beta$ respectively.

We observe that any cube is contained in a dyadic cube $Q_\beta\in \mathcal D^\beta$ for some $\beta\in \{0,\frac{1}{3}\}^n$ with $l(Q_\beta)\le 6l(Q)$ (see for example \cite{PottReg} for the case $n=1$). Thus
$$\frac{1}{|Q|^{1-\frac{\alpha}{nm}}}\int_Q|f|\le \frac{6^{n-\frac{\alpha}{m}}}{|Q_\beta|^{1-\frac{\alpha}{nm}}}\int_{Q_\beta}|f|$$
and consequently,
\Be
\mathcal M_\alpha f\le 6^{nm-\alpha}\sum_{\beta\in \{0,\frac{1}{3}\}^n}\mathcal M_{d,\alpha}^{\beta}f.
\Ee
We will use the following notations: for $\vec {\sigma}=(\sigma_1,\cdots,\sigma_m)$ and 
$\vec f=(f_1,\cdots,f_m)$ given, we write 
$\vec {\sigma}\vec {f}:=(\sigma_1f_1,\cdots,\sigma_mf_m)$ and for a real number $\sigma$, $\sigma \vec {f}=(\sigma f_1,\cdots,\sigma f_m)$. We also use the notation $\vec {\chi_Q}=(\chi_Q,\cdots,\chi_Q)$ ($m$-entries vector) so that
 $\sigma\vec {\chi_Q}=(\sigma\chi_Q,\cdots,\sigma\chi_Q)$ and 
 $\vec {\sigma}\vec {\chi_Q}=(\sigma_1\chi_Q,\cdots,\sigma_m\chi_Q)$.
\subsection{Proof of Proposition \ref{prop:suffsawyertype}}
From the above observations, we see that for the proof of Proposition \ref{prop:suffsawyertype}, it is enough to prove the following.
\bprop\label{prop:dyaSawyer}
Given $\sigma_1,\cdots,\sigma_m$ and $\omega$, $m+1$ weights on $\mathbb R^n$, and $1< p_1,\cdots,p_m<\infty$, let $\frac{1}{p}=\frac{1}{p_1}+\cdots+\frac{1}{p_m}$, $p\le q<\infty$, and define $\nu_{\vec {\sigma}}=\prod_{i=1}^m\sigma_i^{p/p_i}$. Then
if there exists a constant $C_1>0$ such that for any cube $Q\in \mathcal D$, $$\int_Q\left(\mathcal M_{d,\alpha}(\vec {\sigma} \vec {\chi_Q})(x)\right)^q\omega(x)dx\le C_1\left(\nu_{\vec {\sigma}}(Q)\right)^{q/p},$$
then
there exists a constant $C_2>0$ such that $$\int_{\mathbb R^n}\left(\mathcal M_{d,\alpha}(\vec {\sigma} \vec {f})(x)\right)^q\omega(x)dx\le C_2\left(\prod_{i=1}^m\|f\|_{p_i,\sigma_i}\right)^q.$$

Moreover, if $$[\nu_{\vec {\sigma}},\omega]_{\mathcal S_{\vec P,q}}:=\sup_{Q\in \mathcal D}\left(\frac{\int_Q\left(\mathcal M_{d,\alpha}(\vec {\sigma} \vec {\chi_Q})(x)\right)^q\omega(x)dx}{\left(\nu_{\vec {\sigma}}(Q)\right)^{q/p}}\right)^{1/q},$$ then
\Be\label{eq:Spqestim}\|\mathcal M_{d,\alpha}(\vec {\sigma}\cdot)\|_{\left(\prod_{i=1}^mL^{p_i}(\sigma_i)\right)\rightarrow L^q(\omega)}\lesssim [\nu_{\vec {\sigma}},\omega]_{\mathcal S_{\vec P,q}}.\Ee
\eprop

\begin{proof} 
Let $a>2^{nm-\alpha}$. To each integer $k$, we associate the following set $$\Omega_k:=\{x\in \mathbb R^n:a^k<\mathcal M_{d,\alpha}(\sigma\vec {f})(x)\le a^{k+1}\}.$$
There exists a family $\{Q_{k,j}\}_{j\in \mathbb N_0}$ of dyadic cubes maximal with respect to the inclusion and such that $$\prod_{i=1}^m\frac{1}{|Q_{k,j}|^{1-\frac{\alpha}{nm}}}\int_{Q_{k,j}}|f_i(x)|\sigma(x)dx>a^k$$
so that $$\Omega_k\subseteq \cup_{j\in \mathbb N_0}Q_{k,j}.$$
Note that because of their maximality, we have for each fixed $k$, $Q_{k,j}\cap Q_{k,l}=\emptyset$ for $j\neq l$. Also,
$$2^{nm-\alpha}a^k>\prod_{i=1}^m\frac{1}{|Q_{k,j}|^{1-\frac{\alpha}{nm}}}\int_{Q_{k,j}}|f_i(x)|\sigma(x)dx>a^k.$$
Let us put $E(Q_{k,j}):=Q_{k,j}\cap \Omega_k$. Then $\Omega_k=\cup_{j=1}^\infty E_{k,j}$ and the $E(Q_{k,j})$ are disjoint for all $j$ and $k$, i.e
$E(Q_{k,j})\cap E(Q_{l,m})=\emptyset$ for $(k,j)\neq (l,m)$. It follows that
\Beas
\int_{\mathbb R^n}\left(\mathcal M_{d,\alpha}(\vec {\sigma}\vec {f})(x)\right)^q\omega(x)dx &=& \sum_k\int_{\Omega_k}\left(\mathcal M_{d,\alpha}(\vec {\sigma}\vec {f})(x)\right)^q\omega(x)dx\\ &\le&
a^q\sum_ka^{kq}\omega(\Omega_k)\\ &\le& a^q\sum_{k,j}a^{kq}\omega(E(Q_{k,j}))\\ &\le& a^q\sum_{k,j}\left(\prod_{i=1}^m\frac{1}{|Q_{k,j}|^{1-\frac{\alpha}{nm}}}\int_{Q_{k,j}}|f_i(x)|\sigma_i(x)dx\right)^q\omega(E(Q_{k,j}))\\ &=& a^q\sum_{k,j}\omega(E(Q_{k,j}))\left(\prod_{i=1}^m\frac{\sigma_i(Q_{k,j})}{|Q_{k,j}|^{1-\frac{\alpha}{nm}}}\right)^q\left(\prod_{i=1}^mm_{\sigma_i}(f_i,Q_{k,j})\right)^q\\ &=& a^q\sum_{Q\in \mathcal D}\lambda_Q\left(\prod_{i=1}^mm_{\sigma_i}(f_i,Q)\right)^q,
\Eeas
where
$$
\lambda_Q:=\left\{ \begin{matrix} \omega(E(Q))\left(\prod_{i=1}^m\frac{\sigma_i(Q)}{|Q|^{1-\frac{\alpha}{nm}}}\right)^q &\text{if }&\\ Q=Q_{k,j}\,\,\,\textrm{for some}\,\,\,(k,j),\\
      0 & \text{ otherwise} .
                                  \end{matrix} \right.
$$
We observe that for any $R\in \mathcal D$,
\Beas
\sum_{Q\subseteq R, Q\in \mathcal D}\omega(E(Q))\left(\prod_{i=1}^m\frac{\sigma_i(Q)}{|Q|^{1-\frac{\alpha}{nm}}}\right)^q &=& \sum_{Q\subseteq R, Q\in \mathcal D}\omega(E(Q))\left(\prod_{i=1}^m\frac{1}{|Q|^{1-\frac{\alpha}{nm}}}\int_Q\sigma_i\chi_R\right)^q\\ &=& \sum_{Q\subseteq R, Q\in \mathcal D}\int_{E(Q)}\left(\prod_{i=1}^m\frac{1}{|Q|^{1-\frac{\alpha}{nm}}}\int_Q\sigma_i\chi_R\right)^q(x)dx\\ &\le& \int_R\left(\mathcal M_{d,\alpha}(\vec {\sigma}\vec {\chi_R})\right)^q(x)dx\\ &\le&
 [\nu_{\vec {\sigma}},\omega]_{\mathcal S_{\vec P,q}}^q\left(\nu_{\vec {\sigma}}(R)\right)^{q/p}.
\Eeas
That is $\{\lambda_Q\}_{Q\in \mathcal D}$ is a $(\frac{q}{p},\nu_{\vec {\sigma}})$-Carleson sequence. Thus from Lemma \ref{lem:Carlembedchendamian} we obtain
\Beas
\int_{\mathbb R^n}\left(\mathcal M_{d,\alpha}(\vec {\sigma}\vec {f})(x)\right)^q\omega(x)dx &\lesssim& [\nu_{\vec {\sigma}},\omega]_{\mathcal S_{\vec P,q}}^q\prod_{i=1}^m\|M_{d}^{\sigma_i}\left(\sigma_i f_i\right)\|_{p_i,\sigma}^q\\ &\lesssim& [\nu_{\vec {\sigma}},\omega]_{\mathcal S_{\vec P,q}}^q\prod_{i=1}^m\|f_i\|_{p_i,\sigma_i}^q.
\Eeas
The proof is complete.
\end{proof}
\subsection{Proof of Theorem \ref{thm:LiSun}}

We observe again that we only need to prove the following.
\bprop\label{prop:dyaLiSun}
Given $\sigma_1,\cdots,\sigma_m$ and $\omega$, $m+1$ weights on $\mathbb R^n$, and $1< p_1,\cdots,p_m<\infty$, let $\vec {\sigma}=(\sigma_1,\cdots,\sigma_m)$, $\frac{1}{p}=\frac{1}{p_1}+\cdots+\frac{1}{p_m}$ and $q\ge \max\{p_1,\cdots,p_m\}$. Then the following are equivalent.
\begin{itemize}
\item[(i)]
There exists a constant $C_1>0$ such that for any cube $Q\in \mathcal D$, $$\int_Q\left(\mathcal M_{d,\alpha}(\vec {\sigma} \vec {\chi_Q})(x)\right)^q\omega(x)dx\le C_1\prod_{i=1}^m\left(\sigma_i(Q)\right)^{q/{p_i}}.$$
\item[(ii)]
There exists a constant $C_2>0$ such that $$\int_{\mathbb R^n}\left(\mathcal M_{d,\alpha}(\vec {\sigma} \vec {f})(x)\right)^q\omega(x)dx\le C_2\left(\prod_{i=1}^m\|f\|_{p_i,\sigma_i}\right)^q.$$
\end{itemize}
Moreover, if $$[\vec {\sigma},\omega]_{\mathcal S_{\vec P,q}}:=\sup_{Q}\left(\frac{\int_Q\left(\mathcal M_\alpha(\vec {\sigma} \vec {\chi_Q})(x)\right)^q\omega(x)dx}{\prod_{i=1}^m\left(\sigma(Q)\right)^{q/{p_i}}}\right)^{1/q},$$ then
\Be\label{eq:Spqestimlisun}\|\mathcal M_\alpha(\vec {\sigma}\cdot)\|_{\left(\prod_{i=1}^mL^{p_i}(\sigma_i)\right)\rightarrow L^q(\omega)}\backsimeq [\vec {\sigma},\omega]_{\mathcal S_{\vec P,q}}.\Ee
\eprop
As in \cite{LiSun} we restrict ourself to the bilinear case  as the general case follows the same. We will focus on the proof of the sufficiency that is the implication $(\textrm{i})\Rightarrow (\textrm{ii})$ as the converse is obvious. We start by the following lemma proved in \cite{LiSun} and provide a simplified proof.
\blem\label{lem:LiSun1}
Suppose that $0\le \alpha<2n$, that $1<p_1,p_2<\infty$. Put $\frac{1}{p}=\frac{1}{p_1}+\frac{1}{p_2}$ and let $q\ge p_2$, $\sigma_1,\sigma_2,\omega$ be three weights. Then if $f$ is a function with $\supp f\subset R\in \mathcal D$, then
\Be\label{eq:LiSun}
\|\chi_R\mathcal {M}_{d,\alpha}(\chi_R\sigma_1,f\sigma_2)\|_{q,\omega}\lesssim [\vec {\sigma},\omega]_{\mathcal S_{\vec P,q}}\sigma_1(R)^{1/{p_1}}\|f\|_{p_2,\sigma_2}.
\Ee
\elem
\begin{proof}
We proceed as in the proof of Proposition \ref{prop:dyaSawyer}. Let $a>2^{nm-\alpha}$. To each integer $k$, associate the set $$\Omega_k:=\{x\in \mathbb R^n:a^k<\mathcal M_{d,\alpha}(\chi_R\sigma_1,f\sigma_2)(x)\le a^{k+1}\}.$$
There exists a family $\{Q_{k,j}\}_{j\in \mathbb N_0}$ of dyadic cubes maximal with respect to the inclusion and such that  $$\frac{1}{|Q_{k,j}|^{2-\frac{\alpha}{n}}}\int_{Q_{k,j}}\chi_R(x)\sigma_1(x)dx\int_{Q_{k,j}}|f(x)|\sigma_2(x)dx>a^k$$
so that $$\Omega_k\subseteq \cup_{j\in \mathbb N_0}Q_{k,j}.$$
Following the same steps as in the proof of Proposition \ref{prop:dyaSawyer} and using the same notations, we obtain that
\Beas
\int_R\left(\mathcal {M}_{d,\alpha}(\chi_R\sigma_1,f\sigma_2)(x)\right)^q\omega(x)dx &\le& a^q\sum_{Q_{k,j}\subseteq R}\left(\frac{\int_{Q_{k,j}}\chi_R\sigma_1\int_{Q_{k,j}}|f|\sigma_2}{|Q_{k,j}|^{2-\frac{\alpha}{n}}}\right)^q\omega(E(Q_{k,j}))\\ &=& a^q\sum_{Q_{k,j}\subset R}\left(\frac{\int_{Q_{k,j}}\chi_R\sigma_1\int_{Q_{k,j}}|f|\sigma_2}{|Q_{k,j}|^{2-\frac{\alpha}{n}}}\right)^q\omega(E(Q_{k,j}))\\ &+& a^q\left(\frac{\int_{R}\chi_R\sigma_1\int_{R}|f|\sigma_2}{|R|^{2-\frac{\alpha}{n}}}\right)^q\omega(E(R))\\ &=& a^q(T_1+T_2)
\Eeas
where $$T_1=\sum_{Q_{k,j}\subset R}\left(\frac{\int_{Q_{k,j}}\chi_R\sigma_1\int_{Q_{k,j}}|f|\sigma_2}{|Q_{k,j}|^{2-\frac{\alpha}{n}}}\right)^q\omega(E(Q_{k,j})),$$
and
$$T_2=\left(\frac{\int_{R}\chi_R\sigma_1\int_{R}|f|\sigma_2}{|R|^{2-\frac{\alpha}{n}}}\right)^q\omega(E(R)).$$
We easily obtain that
\Beas
T_2 &=& \left(\frac{\int_{R}\chi_R\sigma_1\int_{R}|f|\sigma_2}{|R|^{2-\frac{\alpha}{n}}}\right)^q\omega(E(R))\\ &\le& \left(\frac{1}{\sigma_2(R)}\int_{R}|f|\sigma_2\right)^q\int_R\left(\mathcal {M}_{d,\alpha}(\chi_R\sigma_1,\chi_R\sigma_2)\right)^q\omega(x)dx\\ &\le& [\vec {\sigma},\omega]_{\mathcal S_{\vec P,q}}^q\sigma_1(R)^{q/{p_1}}\|f\|_{p_2,\sigma_2}^q.
\Eeas
We observe that
\Beas
T_1 &=&\sum_{Q_{k,j}\subset R}\left(\frac{\int_{Q_{k,j}}\chi_R\sigma_1\int_{Q_{k,j}}|f|\sigma_2}{|Q_{k,j}|^{2-\frac{\alpha}{n}}}\right)^q\omega(E(Q_{k,j}))\\ &=& \sum_{Q_{k,j}\subset R}m_{\sigma_2}(|f|,Q_{k,j})^q\left(\frac{\int_{Q_{k,j}}\chi_R\sigma_1\int_{Q_{k,j}}\chi_R\sigma_2}
{|Q_{k,j}|^{2-\frac{\alpha}{n}}}\right)^q\omega(E(Q_{k,j}))\\ &=& \sum_{Q\in \mathcal D}\lambda_Qm_{\sigma_2}(|f|,Q)^q
\Eeas
where
$$
\lambda_Q:=\left\{ \begin{matrix} \left(\frac{\int_{Q}\chi_R\sigma_1\int_{Q}\chi_R\sigma_2}
{|Q|^{2-\frac{\alpha}{n}}}\right)^q\omega(E(Q)) &\text{if }&\\ Q=Q_{k,j}\subset R\,\,\,\textrm{for some}\,\,\,(k,j),\\
      0 & \text{ otherwise} .
                                  \end{matrix} \right.
$$
Let us prove that $\{\lambda_Q\}_{Q\in \mathcal D}$ is $(\frac{q}{p_2},\sigma_2)$-Carlseon sequence. Let $K\in \mathcal D$, $K\subset R$. Then
\Beas
\sum_{Q\in \mathcal D, Q\subseteq K}\lambda_Q &=& \sum_{Q\in \mathcal D, Q\subseteq K}\left(\frac{\int_{Q}\chi_R\sigma_1\int_{Q}\chi_R\sigma_2}
{|Q|^{2-\frac{\alpha}{n}}}\right)^q\omega(E(Q))\\ &\le& \sum_{Q\in \mathcal D, Q\subseteq K}\int_{E(Q)}\left(\frac{\int_{Q}\chi_R\sigma_1\int_{Q}\chi_R\sigma_2}
{|Q|^{2-\frac{\alpha}{n}}}\right)^q\omega(x)dx\\ &\le& \int_T\left(\mathcal {M}_{d,\alpha}(\chi_R\sigma_1,\chi_R\sigma_2)(x)\right)^q\omega(x)dx\\ &\lesssim& \sigma_1(T)^{q/{p_1}}\sigma_2(T)^{q/{p_2}}[\vec {\sigma},\omega]_{\mathcal S_{\vec P,q}}^q\\ &\lesssim& \sigma_1(R)^{q/{p_1}}[\vec {\sigma},\omega]_{\mathcal S_{\vec P,q}}^q\sigma_2(T)^{q/{p_2}}.
\Eeas
That is $\{\lambda_Q\}_{Q\in \mathcal D}$ is $(\frac{q}{p_2},\sigma_2)$-Carlseon sequence with constant  $A_{Carl}\lesssim \sigma_1(R)^{q/{p_1}}[\vec {\sigma},\omega]_{\mathcal S_{\vec P,q}}^q$. Thus by Theorem \ref{thm:Carlembed1},
$$T_1\lesssim \sigma_1(R)^{q/{p_1}}[\vec {\sigma},\omega]_{\mathcal S_{\vec P,q}}^q\|f\|_{p_2,\sigma_2}^q.$$
From the estimates of $T_1$ and $T_2$, we conclude that
$$\int_R\left(\mathcal {M}_{d,\alpha}(\chi_R\sigma_1,f\sigma_2)(x)\right)^q\omega(x)dx\lesssim \sigma_1(R)^{q/{p_1}}[\vec {\sigma},\omega]_{\mathcal S_{\vec P,q}}^q\|f\|_{p_2,\sigma_2}^q.$$
The proof is complete.
\end{proof}
The next result is enough to conclude for the sufficient part in Proposition \ref{prop:dyaLiSun} for the bilinear case and was also proved in \cite{LiSun}, we give a simplified proof that uses the general Carleson embedding.
\bprop\label{prop:LiSunsuffdyabili}
Suppose that $0\le \alpha<2n$, that $1<p_1,p_2<\infty$. Put $\frac{1}{p}=\frac{1}{p_1}+\frac{1}{p_2}$ and let $q\ge p_2$, and $\sigma_1,\sigma_2,\omega$ be three weights. Then
\Be\label{eq:LiSunsuffdyabili}
\|\mathcal {M}_{d,\alpha}(f_1\sigma_1,f_2\sigma_2)\|_{q,\omega}\lesssim [\vec {\sigma},\omega]_{\mathcal S_{\vec P,q}}\|f\|_{p_1,\sigma_1}\|f\|_{p_2,\sigma_2}.
\Ee
\eprop
\begin{proof}
From the decomposition in Proposition \ref{prop:dyaSawyer}, we have that
\Beas
\int_{\mathbb R^n}\left(\mathcal {M}_{d,\alpha}(f_1\sigma_1,f_2\sigma_2)(x)\right)^q\omega(x)dx &\le& \sum_{k,j}\left(\frac{\int_{Q_{k,j}}|f_1|\sigma_1\int_{Q_{k,j}}|f_2|\sigma_2}{|Q_{k,j}|^{2-\frac{\alpha}{n}}}\right)^q\omega(E(Q_{k,j}))\\ &=& \sum_{k,j}m_{\sigma_1}(|f_1|,Q_{k,j})^q\left(\frac{\sigma_1(Q_{k,j})\int_{Q_{k,j}}|f_2|\sigma_2}{|Q_{k,j}|^{2-\frac{\alpha}{n}}}\right)^q
\omega(E(Q_{k,j}))\\ &=& \sum_{Q\in \mathcal D}\lambda_Qm_{\sigma_1}(|f_1|,Q)^q
\Eeas
where
$$
\lambda_Q:=\left\{ \begin{matrix} \left(\frac{\sigma_1(Q)\int_{Q}|f_2|\sigma_2}
{|Q|^{2-\frac{\alpha}{n}}}\right)^q\omega(E(Q)) &\text{if }& Q=Q_{k,j}\,\,\,\textrm{for some}\,\,\,(k,j),\\
      0 & \text{ otherwise} .
                                  \end{matrix} \right.
$$
Let us prove that $\{\lambda_Q\}_{Q\in \mathcal D}$ is $(\frac{q}{p_1},\sigma_1)$-Carlseon sequence.

For  $R\in \mathcal D$ given, we obtain using Lemma \ref{lem:LiSun1} that
\Beas
\sum_{Q\in \mathcal D, Q\subseteq R}\lambda_Q &=& \sum_{Q\in \mathcal D, Q\subseteq R}\left(\frac{\sigma_1(Q)\int_{Q}|f_2|\sigma_2}
{|Q|^{2-\frac{\alpha}{n}}}\right)^q\omega(E(Q))\\ &=& \sum_{Q\in \mathcal D, Q\subseteq R}\left(\frac{\sigma_1(Q)\int_{Q}\chi_R|f_2|\sigma_2}
{|Q|^{2-\frac{\alpha}{n}}}\right)^q\omega(E(Q))\\ &\le& \int_R\left(\mathcal {M}_{d,\alpha}(\chi_R\sigma_1,\chi_R|f_2|\sigma_2)\right)^q\omega(x)dx\\ &\le& [\vec {\sigma},\omega]_{\mathcal S_{\vec P,q}}^q\sigma_1(R)^{q/{p_1}}\|f_2\|_{p_2,\sigma_2}^q.
\Eeas
That is $\{\lambda_Q\}_{Q\in \mathcal D}$ is $(\frac{q}{p_1},\sigma_1)$-Carlseon sequence with constant $A_{Carl}\lesssim [\vec {\sigma},\omega]_{\mathcal S_{\vec P,q}}^q\|f_2\|_{p_2,\sigma_2}^q$. Thus using the equivalent definitions in Theorem \ref{thm:Carlembed1}, we obtain
$$\int_{\mathbb R^n}\left(\mathcal {M}_{d,\alpha}(f_1\sigma_1,f_2\sigma_2)(x)\right)^q\omega(x)dx\lesssim [\vec {\sigma},\omega]_{\mathcal S_{\vec P,q}}^q\|f_1\|_{p_1,\sigma_1}^q\|f_2\|_{p_2,\sigma_2}^q.$$
The proof is complete.

\end{proof}

\section{Two-weight norm estimates}
For the proof of Theorem \ref{thm:main}, we also only have to prove the following.
\btheo\label{thm:dyamain}
Let $1<p_1,\cdots,p_m, q<\infty$, and $\vec {\sigma}=(\sigma_1,\cdots,\sigma_m)$, $\omega$ be weights. Put
$\frac{1}{p}=\frac{1}{p_1}+\cdots+\frac{1}{p_m}$ and $\nu_{\vec {\sigma}}=\prod_{i=1}^m\sigma_i^{p/p_i}$. Assume also that $p\le q$.
Then

\begin{align}\label{eq:dyamain1}
\|\mathcal M_\alpha(\vec {\sigma}\vec f)\|_{q,\omega} &\le C(n,p,q)\left([\vec {\sigma},\omega]_{B_{\vec P,q}}\right)^{1/p}\|M_{d}^{\vec {\sigma}}\vec {f}\|_{p,\nu_{\vec \sigma}}\\ &\le C(n,p,q)\left([\vec {\sigma},\omega]_{B_{\vec P,q}}\right)^{1/p}\prod_{i=1}^m\|M_{d}^{\sigma_i}f_i\|_{p_i,\sigma_i}\notag\\ &\le C(n,p,q)\left([\vec {\sigma},\omega]_{B_{\vec P,q}}\right)^{1/p}\prod_{i=1}^m\|f_i\|_{p_i,\sigma_i},\notag
\end{align}

\begin{align}\label{eq:dyamain2}
\|\mathcal M_\alpha(\vec {\sigma}\vec f)\|_{q,\omega} &\le C(n,p,q)\left([\vec {\sigma},\omega]_{A_{\vec P,q}}\right)^{1/p}\left(\prod_{i=1}^m[\sigma_i]_{A_\infty}^{1/p_i}\right)
\prod_{i=1}^m\|M_{d}^{\sigma_i}f_i\|_{p_i,\sigma_i}\\ &\le C(n,p,q)\left([\vec {\sigma},\omega]_{A_{\vec P,q}}\right)^{1/p}\left(\prod_{i=1}^m[\sigma_i]_{A_\infty}^{1/p_i}\right)
\prod_{i=1}^m\|f_i\|_{p_i,\sigma_i}.\notag
\end{align}

\begin{align}\label{eq:dyamain3}
\|\mathcal M_\alpha(\vec {\sigma}\vec f)\|_{q,\omega} &\le C(n,p,q)\left([\vec {\sigma},\omega]_{A_{\vec P,q}}[\vec \sigma]_{W_{\vec P}^\infty}\right)^{1/p}\|M_{d}^{\vec {\sigma}}\vec {f}\|_{p,\nu_{\vec \sigma}}\\ &\le C(n,p,q)\left([\vec {\sigma},\omega]_{A_{\vec P,q}}[\vec \sigma]_{W_{\vec P}^\infty}\right)^{1/p}
\prod_{i=1}^m\|M_{d}^{\sigma_i}f_i\|_{p_i,\sigma_i}\notag\\ &\le C(n,p,q)\left([\vec {\sigma},\omega]_{A_{\vec P,q}}[\vec \sigma]_{W_{\vec P}^\infty}\right)^{1/p}\prod_{i=1}^m\|f_i\|_{p_i,\sigma_i}.\notag
\end{align}
\etheo
Note that the second and the third inequalities in (\ref{eq:dyamain1}), (\ref{eq:dyamain2}) and (\ref{eq:dyamain3}) follow from H\"older's inequality and (\ref{eq:dyamaxfunctestim}) respectively.
\begin{proof}
The proof of Theorem \ref{thm:dyamain} follows essentially as in the case $\alpha=0$ and $p=q$ in \cite{ChenDamian, DamLerPer}. Hence we will only prove (\ref{eq:dyamain2}).

We keep the notations of the proof of Proposition \ref{prop:dyaSawyer}. Let us put $q_i=\frac{qp_i}{p}>1$ and observe that $\frac{1}{q}=\frac{1}{q_1}+\cdots+\frac{1}{q_m}$.

\begin{proof}[{\bf \Proof (\ref{eq:dyamain2})}]. 
We use for simplicity, the notation $Q_{k,j}=Q$ and obtain that
\Beas \int_{\mathbb R^n}\left(\mathcal M_{d,\alpha}(\vec {\sigma}\vec f)(x)\right)^q\omega(x)dx &\lesssim& a^q\sum_{k,j}\left(\prod_{i=1}^m \frac{1}{|Q|^{1-\frac{\alpha}{nm}}}\int_Q|f_i|\sigma_i\right)^q\omega(Q)\\ &=& \sum_{k,j}\omega(Q)\left(\prod_{i=1}^m\frac{\sigma_i(Q)}{|Q|^{1-\frac{\alpha}{nm}}}\right)^q
\prod_{i=1}^m\left(m_{\sigma_i}(|f_i|,Q)\right)^q\\ &=& \sum_{k,j }\left(\frac{\omega(Q)^{p/q}\prod_{i=1}^m\sigma_i(Q)^{p/{p_i'}}}{|Q|^{p(m-\frac{\alpha}{n})}}\right)^{q/p}
\prod_{i=1}^m\sigma_i(Q)^{q/p_i}\left(m_{\sigma_i}(|f_i|,Q)\right)^q\\ &\le& [\vec {\sigma},\omega]_{A_{\vec P,q}}^{q/p}\sum_{k,j}\prod_{i=1}^m\sigma_i(Q)^{q/p_i}\left(m_{\sigma_i}(|f_i|,Q)\right)^q\\ &\le& [\vec {\sigma},\omega]_{A_{\vec P,q}}^{q/p}\prod_{i=1}^m\left(\sum_{k,j}\sigma_i(Q)^{q_i/p_i}\left(m_{\sigma_i}(|f_i|,Q)\right)^{q_i}\right)^{q/q_i}\\ &=& [\vec {\sigma},\omega]_{A_{\vec P,q}}^{q/p}\prod_{i=1}^m\left(\sum_{k,j}\sigma_i(Q)\left(m_{\sigma_i}(|f_i|,Q)\right)^{p_i}\right)^{q/p_i}\\ &=& [\vec {\sigma},\omega]_{A_{\vec P,q}}^{q/p}\prod_{i=1}^m\left(\sum_{Q\in \mathcal D}\lambda_Q^i\left(m_{\sigma_i}(|f_i|,Q)\right)^{p_i}\right)^{q/p_i}
\Eeas
with
$$
\lambda_Q^i:=\left\{ \begin{matrix} \sigma_i(Q) &\text{if }& Q=Q_{k,j}\,\,\,\textrm{for some}\,\,\,(k,j),\\
      0 & \text{ otherwise} .
                                  \end{matrix} \right.
$$
By Lemma \ref{lem:Carlembedlemma} we can conclude that
\Beas L &:=& \int_{\mathbb R^n}\left(\mathcal M_{d,\alpha}(\vec {\sigma}\vec f)(x)\right)^q\omega(x)dx\\ &\le& C(n,p,q)\left([\vec {\sigma},\omega]_{A_{\vec P,q}}\right)^{q/p}\left(\prod_{i=1}^m[\sigma_i]_{A_\infty}^{q/p_i}\right)
\prod_{i=1}^m\|M_{d}^{\sigma_i}f_i\|_{p_i,\sigma_i}^q
\Eeas
provided that for each $i=1,\cdots,m$,$\{\lambda_Q^i\}_{Q\in \mathcal D}$ is a $\sigma_i$-Carleson sequence with the appropriate constant. That is for any $R\in \mathcal D$,
\Be\label{eq:Carlembeddyasawyer}
\sum_{Q\subseteq R, Q\in \mathcal D}\sigma_i(Q)\lesssim [\sigma_i]_{A_\infty}\sigma(R).
\Ee

The inequality (\ref{eq:Carlembeddyasawyer}) can be found in \cite{HytPerez}. Let us prove it here for completeness.  

We first check the following.
\Be\label{eq:QEQestim}
|Q_{k,j}|\le \gamma |E(Q_{k,j})|,\,\,\,\gamma=\frac{1}{1-\left(\frac{2^{nm-\alpha}}{a}\right)^{1/(m-\frac{\alpha}{n})}}.
\Ee
Recall that $E(Q_{k,j})=Q_{k,j}\cap \Omega_k$ and observe that
\Beas
|Q_{k,j}\cap\left(\cup_{l=1}^\infty Q_{k+1,l}\right)| &=& \sum_{Q_{k+1,l}\subset Q_{k,j}}|Q_{k+1,l}|\\ &\le& \left(\sum_{Q_{k+1,l}\subset Q_{k,j}}|Q_{k+1,l}|^{m-\frac{\alpha}{n}}\right)^{\frac{1}{m-\frac{\alpha}{n}}}\\ &\le& \left(\frac{1}{a^{k+1}}\sum_{Q_{k+1,l}\subset Q_{k,j}}\prod_{i=1}^m\int_{Q_{k+1,l}}|f_i|\sigma_i\right)^{\frac{1}{m-\frac{\alpha}{n}}}\\ &\le& \left(\frac{1}{a^{k+1}}\prod_{i=1}^m\sum_{Q_{k+1,l}\subset Q_{k,j}}\int_{Q_{k+1,l}}|f_i|\sigma_i\right)^{\frac{1}{m-\frac{\alpha}{n}}}\\ &\le&  \left(\frac{1}{a^{k+1}}\prod_{i=1}^m\int_{Q_{k,j}}|f_i|\sigma_i\right)^{\frac{1}{m-\frac{\alpha}{n}}}\\ &\le& \left(\frac{1}{a^{k+1}}2^{nm-\alpha}a^k|Q_{k,j}|^{m-\frac{\alpha}{n}}\right)^{\frac{1}{m-\frac{\alpha}{n}}}\\ &\le& \left(\frac{2^{nm-\alpha}}{a}\right)^{1/(m-\frac{\alpha}{n})}|Q_{k,j}|.
\Eeas
Thus
\Beas
|Q_{k,j}| &\le& E(Q_{k,j})+|Q_{k,j}\cap\left(\cup_{l=1}^\infty Q_{k+1,l}\right)|\\ &\le& |E(Q_{k,j})|+ \left(\frac{2^{nm-\alpha}}{a}\right)^{1/(m-\frac{\alpha}{n})}|Q_{k,j}|
\Eeas
which proves (\ref{eq:QEQestim}). It follows that
\Beas
\sum_{Q\subseteq R, Q\in \mathcal D}\sigma(Q) &\le& \gamma\sum_{Q\subseteq R, Q\in \mathcal D}\frac{\sigma_i(Q)}{|Q|}|E(Q)|\\ &=& \gamma\sum_{Q\subseteq R, Q\in \mathcal D}\int_{E(Q)}\frac{\sigma_i(Q)}{|Q|}\\ &\le& \gamma\sum_{Q\subseteq R, Q\in \mathcal D}\int_{E(Q)}M_d(\sigma_i\chi_R)\\ &=& \gamma\int_RM_d(\sigma_i\chi_R)\\ &\le& \gamma[\sigma_i]_{A_\infty}\sigma_i(R).
\Eeas
\end{proof}
The proof is complete.
\end{proof}
\bibliographystyle{plain}

\end{document}